\documentclass{birkjour}
\usepackage[all]{xy}
\usepackage{wrapfig}
\usepackage{url}

\usepackage{epic}
\usepackage{amsmath}[1996/11/01]
\usepackage{amssymb,amsthm,amsfonts,latexsym,epsfig,graphics}
 \def\beql#1#2\eeql{\begin{equation}\label{#1}#2\end{equation}}

\DeclareMathOperator{\Mat}{Mat}
\DeclareMathOperator{\Double}{Double}
\DeclareMathOperator{\Inv}{Inv}

\DeclareMathOperator{\Gal}{Gal}
\DeclareMathOperator{\Bil}{Bil}

\DeclareMathOperator{\Aut}{Aut}

\DeclareMathOperator{\ccwe}{ccwe}
\DeclareMathOperator{\fwe}{fwe}

\DeclareMathOperator{\End}{End}
\DeclareMathOperator{\Char}{char}

\DeclareMathOperator{\GL}{GL}

\DeclareMathOperator{\Hom}{Hom}
\DeclareMathOperator{\diag}{diag}

\DeclareMathOperator{\id}{id}

\newcommand{\leftBra}{\mbox{\,$\{ \hspace*{-.09in} \{ $\,}}
\newcommand{\rightBra}{\mbox{\,$\} \hspace*{-.09in} \} $\,}}

\theoremstyle{plain}
\newtheorem{theorem}{Theorem}
\newtheorem{lemma}[theorem]{Lemma}

\newtheorem{corollary}[theorem]{Corollary}
\newtheorem{definition}[theorem]{Definition}
\theoremstyle{remark}

\newtheorem{remark}[theorem]{Remark}

\numberwithin{theorem}{section}

\newcommand{\Z}{{\mathbb{Z}}}
\newcommand{\Q}{{\mathbb{Q}}}
\newcommand{\F}{{\mathbb{F}}}

\newcommand{\N}{{\mathbb{N}}}
\newcommand{\R}{{\mathbb{R}}}

\newcommand{\C}{{\mathbb{C}}}

\newcommand{\II}{{\textsc{II}}} 
\renewcommand{\em}{\sf}

\begin{document}
 \bibliographystyle{plain}

\title{$\Gamma $-conjugate weight enumerators and invariant theory}
\author{Gabriele Nebe and Leonie Scheeren}
\email{nebe@math.rwth-aachen.de, leonie.scheeren@rwth-aachen.de}
\address{Lehrstuhl f\"ur Algebra und Zahlentheorie, RWTH Aachen University, 52056 Aachen, Germany}

\begin{abstract}
	Let $K$ be a field, $\Gamma $ a finite group of 
	field automorphisms of $K$, $F$ the $\Gamma $-fixed field 
	in $K$ and $G\leq \GL_v(K)$ a finite matrix group. 
Then the action of $\Gamma $ defines a grading on the symmetric algebra of 
the $F$-space $K^v$ which we use to introduce 
the notion of homogeneous $\Gamma $-conjugate invariants of $G$. 
We apply this new grading in invariant theory to 
broaden the connection between codes
	and invariant theory by introducing $\Gamma $-conjugate complete 
	weight enumerators of codes. 
The main result of this paper applies the theory from Nebe, Rains, Sloane 
to show that under certain extra conditions these new weight enumerators 
generate the ring of $\Gamma $-conjugate invariants of the associated 
Clifford-Weil groups. As an immediate consequence we 
	obtain a  result by Bannai etal that the 
complex conjugate weight enumerators generate the ring of complex conjugate
invariants of the complex Clifford group.
 Also the Schur-Weyl duality conjectured and 
partly shown by Gross etal can be derived from our main result. 
\\
	{\sc keywords:}  Galois action, invariant ring, 
	generalised Molien series, Clifford-Weil group, self-dual codes, 
	complex projective design, Schur-Weyl duality.
\end{abstract}

\maketitle

\section{Introduction}
\label{sec:intro}

 A complex conjugate polynomial of degree $(N,N)$ in $v$ variables 
 is a homogeneous
 polynomial $p \in \C[x_1,\ldots , x_v, \overline{x}_1,\ldots , \overline{x}_v]$ that is of degree $N$ in the variables $x_1,\ldots , x_v$ and of 
 degree $N$ in their complex conjugates. 
 An invariant theory for
 complex conjugate polynomials has been developed in \cite{Forger}. 
 Given a finite complex unitary matrix group $G\leq \mathrm{U}_v(\C) $ and some 
 $t\in \N$ such that for all $N=1,\ldots ,t$ all complex conjugate 
 invariants of degree $(N,N)$ of $G$ are multiples of the $N$th
 power of the invariant Hermitian form, then all $G$-orbits on $\C^v$ 
 define projective $t$-designs. 
 More generally if $x\in \C^v$ is a zero of all harmonic 
 $G$-invariant polynomials up to degree $(t,t)$, then the $G$-orbit $xG$ 
gives rise to a projective $t$-design. 

The textbook \cite{NRS} gives a very general notion of a 
Type of a code. Associated to a Type $\rho $ and an integer $m\geq 1$ 
 there is a 
finite complex matrix group ${\mathcal C}_m(\rho )$, the 
associated Clifford-Weil group of genus $m$, 
such that the genus-$m$ complete weight enumerator of 
any code of Type $\rho $ and length $N$ is an invariant 
polynomial of ${\mathcal C}_m(\rho )$, homogeneous of degree $N$.
The Weight Enumerator Conjecture states that the 
space of homogeneous degree $N$ invariants of ${\mathcal C}_m(\rho )$ 
is spanned by these weight enumerators. 
If $v$ is the size of the alphabet of the codes of Type $\rho $, then 
${\mathcal C}_m(\rho )$ consists of matrices of size $v^m$. 
Despite of this exponentially growing dimension, the space of 
invariants of a given degree $N$ can be obtained by enumerating 
all codes of length $N$ of Type  $\rho $.

For the Type of doubly even binary codes, where $v=2$, 
the associated Clifford-Weil groups are the complex 
Clifford groups 
 ${\mathcal X}_m \leq \GL_{2^m}(\C) $ which have a tight 
 connection to quantum information theory (see \cite{Zhu}). 
 The main result of  \cite{Bannai} shows that the ring of 
 complex conjugate invariants of ${\mathcal X}_m$
 is spanned 
 by the genus-$m$ complex conjugate weight enumerators of self-dual doubly even binary codes. The paper also enumerates all such codes up to 
 length $(5,5)$ therewith proving a conjecture from \cite{Zhu} that
 whenever an ${\mathcal X}_m$-orbit forms a projective $4$-design 
 then it is automatically a projective $5$-design. 

 We extend the approach in
 \cite{Bannai} to the more general set-up in \cite{NRS}. 
Any Type $\rho $ also determines an abelian number field $K$ 
such that the associated Clifford-Weil groups consist of 
matrices over $K$. 
Given a subgroup $\Gamma $ of the automorphism group of $K$, 
we introduce the concept of $\Gamma $-conjugate weight enumerators of codes of 
 Type $\rho$. These are elements of the ring of $\Gamma $-conjugate invariants 
(Definition \ref{conjpol}) of ${\mathcal C}_m(\rho )$. 
Our main result, Theorem \ref{main}, shows that the space of 
$\Gamma $-conjugate homogeneous invariants of ${\mathcal C}_m(\rho )$ 
is spanned by  these weight enumerators of codes of a given length. 
Therefore the enumeration of such codes of small length 
determines the $\Gamma $-conjugate invariants of small degree 
of the genus-$m$ Clifford-Weil group
associated to  $\rho $ for all $m\in \N$. 

With a view to possible applications to
measurement schemes for low rank matrix recovery from complex projective $t$-designs as in \cite{Rauhut},
Section \ref{examples} gives a few examples of self-dual codes of Type 
$(N,N)$ for $N\leq 4$ for several small representations $\rho $. 
Section \ref{SchurWeyl} contains, as another application of Theorem \ref{main},
a short proof of the Schur-Weyl duality conjectured and partly established 
in \cite{Gross}.

The authors acknowledge funding by 
the Deutsche Forschungsgemeinschaft (DFG, German Research Foundation)
Project number 442047500 through the Collaborative Research Center
``Sparsity and Singular Structures'' (SFB 1481).

 \section{Conjugate polynomials} \label{Secconjpol}

 
 Let $K$ be a field
  and let $\Gamma := \{ \alpha _1,\ldots , \alpha _n \} \leq \Aut(K)$ 
be a group of 
field automorphisms of $K$ that is finite of order $n$.
 Put 
 $$F:= \{ a\in K \mid \alpha _i(a) = a \mbox{ for all } 
 i= 1,\ldots , n \}  $$
 to denote the fixed field of $\Gamma $.
 Then $K/F$ is a Galois extension with Galois group $\Gamma $. 

 For a  vector space $U$ over $K$ of finite dimension $v$ 
 we denote by $K[U]$ the 
 ring of polynomials on $U$. 
 We now restrict scalars and write $U_F$ for $U$ regarded
as an $nv$-dimensional space over $F$. 
 For  a $K$-basis $(u_1,\ldots , u_v)$ 
 of $U$ and an $F$-basis $(b_1,\ldots ,b_n)$ of $K$ 
 the tuple 
 $B:=(b_ju_i : 1\leq j \leq n, 1\leq i \leq v)$ is an $F$-basis of $U_F$.
 This yields an $F$-algebra isomorphism between
the symmetric algebra $F[U_F]$ of the
 dual space $U_F^*$ and $F[y_{ji} :  1\leq j \leq n, 1\leq i \leq v ]$,
 where the functions $y_{ji} : U_F \to F$ form the dual basis 
 to the previously chosen basis $B$. 
The latter is a polynomial ring in $nv$ variables over $F$ 
 with well studied gradings given by multi-degrees.
 As we are mostly interested in invariant theory of finite groups
 over fields of characteristic 0 
 we assume that $K$ is an infinite field and see elements of 
 $K[U] \cong K[x_1,\ldots , x_v]$ as polynomial functions on $U\cong K^v$.

 \begin{definition}\label{conjpol}
	 For $v\in \N$ we denote by $\underline{v}:=\{1,\ldots, v\}$ 
	 and by $K^{K^v}$ the $K$-algebra of $K$-valued functions on $K^v$.
	 Assume that $K$ is an infinite field.
 The {\em ring of $\Gamma $-conjugate polynomials over $K$ 
 in $v$ variables}
 ${\bf x} := (x_1,\ldots , x_v)$ is denoted by
 $$K[ {\bf x} \circ \Gamma ] := 
 K[ x_i \circ \alpha _j : 1\leq i \leq v, 
	 1\leq j \leq n ] \leq K^{K^v} $$
	 and defined as
 the $K$-subspace 
spanned by the monomial functions 
	 $M((m_{ij} \mid (i,j) \in \underline{v}\times \underline{n}) )
	 \in K^{K^v}$, where
	 $$M((m_{ij} \mid (i,j) \in \underline{v}\times \underline{n}) ):= \prod _{j=1}^n \prod _{i=1}^v (x_{i}\circ \alpha_j )^{m_{ij}} \colon 
 (k_1,\ldots , k_v) \mapsto 
 \prod _{j=1}^n \prod _{i=1}^v (\alpha_j (k_i) )^{m_{ij}} .$$
The degree of such a  monomial is 
$$ \deg( M((m_{ij} \mid (i,j) \in \underline{v}\times \underline{n}) ) := (d_1,\ldots , d_n) $$ 
	 with $d_j = \sum _{i=1}^v m_{ij} $ for all $j\in \underline{n}$. 
 \end{definition} 

The $K$-algebra structure of $K[{\bf x}\circ \Gamma ]$ 
is inherited from the $K$-algebra structure of $K^{K^v}$.
In particular the multiplication of two monomials is given by 
 $$M((\ell_{ij} \mid (i,j) \in \underline{v}\times \underline{n})  ) 
 M((m_{ij} \mid (i,j) \in \underline{v}\times \underline{n}) )  
 = M ((\ell_{ij} + m_{ij}\mid (i,j) \in \underline{v}\times \underline{n}) )  $$
and as usual the degree of the product is just the sum of the degrees 
of the two factors. 

The $\Gamma $-action 
provides a finer notion of degree of a polynomial in $K[U_F] := K\otimes _{F} F[U_F]$. 
Note that also
  the functions $(x_i \circ \alpha _j \mid (i,j) \in \underline{v} \times \underline{n} )$
 form a basis of the dual space $K\otimes _{F} U_F^*$.
 So we obtain a
  $K$-algebra isomorphism
  \begin{equation} \label{phi}
\varphi : K[{\bf x} \circ \Gamma ]  \to K[y_{ji}  \mid (i,j) \in \underline{v}\times \underline{n} ].
  \end{equation}

This shows that the linear functions $(x_i\circ \alpha _j : (i,j) \in \underline{v} \times \underline{n} ) $ in $K^{K^v}$ are algebraically independent over $K$.



\begin{corollary}
The space of homogeneous polynomials of degree $(d_1,\ldots ,d_n)$ 
in $K[ {\bf x} \circ \Gamma ]$
is the span of all monomials 
	$M((m_{ij} \mid (i,j) \in \underline{v} \times \underline{n}))$ with $d_j = \sum _{i=1}^v m_{ij} $. 
As these monomials form a basis its dimension is 
	$$\dim (K[{\mathbf x}\circ \Gamma ] _{d_1,\ldots , d_n}) = 
	\prod _{j=1}^n \dim (K[x_1,\ldots , x_v]_{d_j}) = 
	\prod _{j=1}^n {{d_j+v-1}\choose{d_j}} .$$
\end{corollary}

\section{Invariant Theory}

We keep the assumptions of the previous section, in particular 
$K$ is an infinite field and $\Gamma =\{\alpha _1,\ldots ,\alpha _n\}$ is
a group of automorphisms 
of $K$ of finite order $n$. 
Let $G\leq \GL_v(K)$ be a group. Then the right action of 
$G$ on $K^v$ defines a  right action of $G$ by $K$-algebra automorphisms
 on the $K$-algebra of $K$-valued functions on ${K^v}$ by
 $$ f\cdot  g: K^v \to K, k \mapsto f(kg^{-1}) \mbox{ for all } 
 g \in G, f\in K^{K^v} .$$
This action preserves the subalgebra 
$K[ {\bf x} \circ \Gamma ] $
as well as its subspaces of homogeneous polynomials
 of a given degree. 

\begin{definition}
Let $K[ {\bf x} \circ \Gamma ]^G $ denote the $K$-algebra of 
$G$-invariant functions in $K[ {\bf x} \circ \Gamma ] $. 
\end{definition}

Via the isomorphism $\varphi $ from Equation \eqref{phi} 
the $K$-algebra of $G$-invariant functions 
$K[ {\bf x} \circ \Gamma ]^G  $ is isomorphic to the 
ring of $G$-invariant polynomials in 
$K[y_{ji}  \mid (i,j) \in \underline{v} \times \underline{n}] $. 
In particular 
 classical invariant theory gives us Molien's formula 
 for the Hilbert series of this invariant ring. 

The grading from Definition \ref{conjpol} 
refines the classical degree function thus giving a 
notion of $\Gamma $-conjugate Hilbert series of the invariant ring 
$K[ {\bf x} \circ \Gamma ]^G  $, generalising the 
Forger series (where $K=\C $ and $F=\R $) \cite{Forger} 
to our situation.

\begin{definition}
	For any ${\bf d}:=(d_1,\ldots , d_n) \in \N _{0}^n$ let
	$$a_{{\bf d}} := \dim (K[{\bf x} \circ \Gamma ]_{{\bf d}} \cap K[ {\bf x} \circ \Gamma ]^G ) $$ 
	denote the dimension of the space of $G$-invariant conjugate
	polynomials that are homogeneous of degree ${\bf d}$. 
	Putting ${\bf z}^{\bf d} :=  z_1^{d_1} \ldots z_n ^{d_n} $ we define
	$${\mathcal H}(K[ {\bf x} \circ \Gamma ]^G  ) :=
	\sum _{{\bf d} \in  \N _{0}^n} a_{{\bf d}} {\bf z}^{{\bf d}} \in \Z [[ z_1,\ldots , z_n ]] $$ 
	 the 
	{\em $\Gamma $-conjugate Hilbert series} 
	of the ring $ K[ {\bf x} \circ \Gamma ]^G $.
\end{definition}

\subsection{Molien's theorem}

If $G$ is finite and $K$ has characteristic $0$, then Molien's theorem 
\cite[Theorem 2.5.2]{Benson} gives a useful expression of 
the Hilbert series of the classical invariant ring of $G$. 
With a completely analogous proof (see also \cite{Forger} for $K=\C $ and 
$F=\R $) we obtain the following theorem.

\begin{theorem}
	Let $K$ be a field of characteristic 0 and let
	$G$ be a finite subgroup of $\GL_v(K)$.
	Then  $${\mathcal H}(K[ {\bf x} \circ \Gamma ]^G  )  = \frac{1}{|G|} 
	\sum_{g\in G}\prod_{j=1}^n \frac{1}{\det(I_v - z_{j}\alpha_j(g))}.$$
\end{theorem}

\begin{proof}
	For a given degree ${\bf d}$ $\in$ $\N _{0}^n$  
	the Reynolds operator 
	$\frac{1}{| G|}\displaystyle\sum_{g \in G}g \in K[G]$
	induces a $K$-linear projection $P_{{\bf d}}$ 
	from $K[{\bf x} \circ \Gamma ]_{{\bf d}}$ onto the 
	fixed space $K[ {\bf x} \circ \Gamma ]_{\bf d}^G $.
	Since char($K$) = 0 the dimension of this fixed subspace 
	is equal to the trace of $P_{{\bf d}}$.  \\
        Therefore 
	$${\mathcal H}(K[ {\bf x} \circ \Gamma ]^G  ) =
	\sum _{{\bf d} \in  \N _{0}^n} \textrm{trace}(P_{{\bf d}}) {\bf z}^{{\bf d}} .$$
	As  $\textrm{trace}(P_{{\bf d}}) = \frac{1}{|G|} \sum _{g\in G} t_{{\bf d}} (g) $ 
it suffices to compute the trace  $t_{{\bf d}} (g)$ 
of the action of $g$ on 
	$K[{\bf x} \circ \Gamma ]_{{\bf d}}$ for all $g\in G$. 
	\\
	To do so we may and will assume 
	that $K$ contains a $|G|$th primitive root of unity. 
   Then each $g$ $\in$ $G$ is diagonalizable over $K$ so after 
	a suitable choice of basis we assume that 
	$g^{-1} = \diag (\lambda _1,\ldots , \lambda _v)$.
	Then for 
	$\alpha_{j}$ $\in \Gamma$, 
	$$\det(I_v - z_j \alpha _j(g)) = \prod_{i = 1}^v (1 - \alpha_j(\lambda_i)z_j).$$ 
	The 
 monomials from Definition \ref{conjpol} 
	form an
eigenvector basis for the action of $g$ on 
	$K[{\bf x} \circ \Gamma ]_{{\bf d}}$, 
	where $M((m_{ij}  \mid (i,j) \in \underline{v} \times \underline{n}))$ has eigenvalue 
	$\displaystyle\prod_{j=1}^{n} \prod_{i=1}^v \alpha_j(\lambda_{i})^{m_{ij}}$.
	So we get 
	$$
		\sum _{{\bf d} \in  \N _{0}^n} t_{{\bf d}}(g) {\bf z}^{{\bf d}} =  \prod_{j=1}^{n} \prod_{i=1}^v 
	\sum_{m\in \N_0}(\alpha_j(\lambda_i)z_{j})^{m} =  $$ 
	$$\prod_{j=1}^n 	\prod_{i=1}^v \frac{1}{1 - \alpha _j(\lambda_{i})z_{j}} = \prod_{j=1}^n \frac{1}{\det(I_v - z_j \alpha_j(g))}.
	 $$
\end{proof}

\section{The Type of a self-dual code}

 This section briefly recalls the relevant notions from \cite{NRS}. 
 Classically a self-dual code $C$  of length $N$ 
 over a finite field $\F _q$ is a linear subspace $C\leq \F_q^N$ 
 that is self-dual (i.e. $C=C^{\perp }$) 
 with respect to the standard inner product. 
 Loosely speaking, to define self-dual codes in a more general sense 
 we need a ring $R$, a left $R$-module $V$ and a 
 non-singular form $\beta $ on $V$ such that $\beta $ can be used 
 to define the orthogonal code $C^{\perp }$ of an $R$-submodule $C\leq V^N$. 
 Extra conditions (such as being doubly even)
 can be imposed by means of  isotropy conditions with respect to 
 some set $Q$ of quadratic maps on $V$. 
 We call such an admissible  quadruple $\rho := (R,V,Q,\beta )$ 
 a {\em Type}, and the self-dual isotropic codes  that 
 arise from $\rho $ are 
 called {\em codes of Type $\rho$}. 
 More precisely, a Type is a representation of an abstract form ring. 

 \subsection{Form rings} 

 A {\em form ring} $(R,M,\psi,\Phi )$ is a quadruple, where 
 $R$ is a (unital and associative) ring, $M$ a right $R\otimes R$-module,
 $\psi : R_R \to M_{1\otimes R}$ an isomorphism of right $R$-modules.
 It comes with an involution $\tau : M \to M$ such that 
 $\tau (m) (s\otimes r) = \tau (m(r\otimes s)) $ for all $m\in M$, $r,s\in R$
 and such that 
 $\epsilon := \psi ^{-1} (\tau (\psi(1))) $ is a unit in $R$. 
 The isomorphism $\psi $ defines an anti automorphism $\ ^J: R \to R$ 
 by $r^J := \psi ^{-1}(\psi(1) (r\otimes 1)) $. 
 The last ingredient is an $R$-qmodule $\Phi $, i.e. an abelian group $\Phi $
 together with a (pointed quadratic) map $[]:  R  \to  \End_{\Z} (\Phi)$ 
 such that $[1]=1$ and $[rs] = [r][s]$ for all $r,s \in R $. 
 There are qmodule homomorphisms $\leftBra \ \rightBra : M \to \Phi $ 
 and $\lambda : \Phi \to M$
 such that for all $m\in M$, $\phi \in \Phi $, 
 $$\leftBra \tau(m) \rightBra = \leftBra m \rightBra ,\ 
 \tau(\lambda(\phi )) = \lambda (\phi ), \ \lambda(\leftBra m \rightBra ) = 
 m + \tau (m)  $$ and 
 $$\leftBra \lambda(\phi ) (r\otimes s) \rightBra = 
 \phi [r+s] - \phi[r] - \phi [s] \mbox{ for all } r,s\in R, \phi \in \Phi .$$
 Taking $M=R$ and $\psi $ the identity 
 we abbreviate $(R,\Phi ):= (R,M,\psi ,\Phi )$.

 \subsection{Finite representations of form rings}

 A finite representation of a form ring $(R, \Phi )$ 
 is a quadruple $\rho = (V,\rho _M, \rho _{\Phi }, \beta )$,
 where $V$  is a 
 left $R$-module of finite cardinality $v$, 
 $\rho _M:M\to \Bil (V)$ is an $R\otimes R$-module 
 homomorphism into the group $\Bil (V)$ of bi-additive $\Q/\Z $-valued 
 maps on $V$ 
 compatible with the involution $\tau $, i.e. 
 $\rho_M(\tau(m))(x,y) = \rho_M(m) (y,x)$ for all 
 $m\in M, x,y\in V$, and such that $\beta := \rho_M(\psi(1)) $ is
 non-singular. 
 Also $\rho _{\Phi }$ is an $R$-qmodule homomorphism from $\Phi $ into the group 
 of $\Q/\Z $-valued quadratic maps on $V$
 satisfying 
 $$\begin{array}{l} 
	 \rho_\Phi (\leftBra m \rightBra ) (x) = \rho _M(m) (x,x) \mbox{ and }
\\ \rho_M(\lambda (\phi )) = \rho _{\Phi} (\phi ) (x+y) - 
  \rho _{\Phi} (\phi ) (x) - 
 \rho _{\Phi} (\phi ) (y)  \end{array} $$
   for all $ x,y\in V, m\in M, \phi \in \Phi .$ 

 \begin{definition}
	 Let $\rho = (V,\rho_M,\rho_{\Phi } , \beta) $ be a finite representation of a form ring $(R,\Phi)$.
	 Any $R$-submodule $C\leq V$ is called a {\em code in $\rho $}.
	 For a code $C$ in $\rho $  the orthogonal module is
	 $$C^{\perp} = C^{\perp , \beta } := 
	 \{ x\in V \mid  \beta (x,c ) = 0 \mbox{ for all } c\in C \} .$$
	 The code $C$ is called {\em self-dual} if $C=C^{\perp} $ and
	 {\em self-orthogonal} if $C\subseteq C^{\perp }$.
	 A self-orthogonal code $C$ in $\rho $ 
	 is called {\em isotropic} if $\rho_{\Phi}(\phi) (C) = \{0 \} $ for all $\phi \in \Phi $. 

	 A {\em code of Type $\rho $} is a self-dual isotropic code in $\rho $. 
 \end{definition}

  MacWilliams transformations (\cite[Section 2.2]{NRS}) 
  map the weight enumerator of a 
  code to the one of its dual and fix weight enumerators
  of self-dual codes. To obtain a generating set of the associated
  Clifford-Weil group we need to include MacWilliams transformations
  for representatives of the conjugacy classes of primitive
  symmetric idempotents.

 \begin{definition}\label{symid}
 An idempotent $\iota ^2=\iota \in R$ is called {\em symmetric}
 if $\iota R \cong \iota ^J R$ as right $R$-modules. 
 Such an isomorphism is given by left multiplication with 
	 some $v_{\iota } \in \iota ^JR \iota$ with inverse
	  $u_{\iota } \in \iota R \iota ^J$ such that 
$u_{\iota } v_{\iota } = \iota $ and $v_{\iota } u_{\iota } = \iota ^J$ 
(see \cite[Section 3.5.3]{NRS}). 
 \end{definition}

 \begin{remark} \label{iota} (see \cite[Theorem 3.5.9]{NRS}) 
	 Let $\iota \in R$ be a symmetric idempotent and $C=C^{\perp} \leq V$
	 a self-dual code. 
	 Then $\iota C$ is a self-dual code in 
	 $\iota V$ and in particular $|\iota C|^2  = |\iota V|$.
 \end{remark}

 \begin{definition}
	 The {\em value group} of the representation $\rho =(V,\rho _M, \rho _{\Phi }, \beta ) $
	 is the subgroup $\nu (\rho )$ of $\Q/\Z $ 
	 generated by 
	 $$ \{ \beta (x,y) \mid x,y \in V \} \cup  
	 \{ \rho _{\Phi }(\phi ) (x) \mid x\in V, \phi \in \Phi \}.$$ 
 As $\nu (\rho )$ is a finitely generated (and hence finite) subgroup of 
	 $\Q/\Z $ it is cyclic, so $\nu (\rho ) = \langle \frac{1}{f} +\Z \rangle $,  where  $f= f (\rho )= |\nu (\rho )|   $ is called the {\em conductor} of $\rho $. 
 \end{definition}

\begin{remark}\label{sum}
	Put $f:=f(\rho)$ and let 
	 $a_1,\ldots , a_N \in \Z $ be prime to $f$ and put 
	 ${\bf a}:= (a_1,\ldots , a_N)$.
	 Then the orthogonal sum $$\rho^{{\bf a}} :=  
			(V^N, a_1\rho_M\perp \ldots \perp a_N \rho_M , a_1\rho_{\Phi} \perp \ldots \perp a_N \rho _{\Phi} , 
			a_1\beta \perp \ldots \perp a_N \beta )$$
			is a finite representation of 
			the form ring $(R,\Phi )$.
\end{remark}

 \section{Clifford-Weil groups and full weight enumerators} \label{CliffordWeil}

 Let 
 $\rho = (V,\rho _M, \rho _{\Phi }, \beta )$
be  a finite representation  of a form ring and put $v:=|V|$. 
The group algebra $\C V$ is a $v$-dimensional complex vector space
with basis $(b_w : w \in V)$. 
The {\em full weight enumerator} of a code in $\rho $ is defined as
$$\fwe(C) := \sum _{c\in C} b_c \in \C V  .$$ 
The associated Clifford-Weil group ${\mathcal C}(\rho )$ 
is a group of linear operators on $\C V$ whose generators are 
explicitly given 
 in \cite[Definition 5.3.1]{NRS}; these generators are 
the obvious transformations that stabilise $\fwe (C)$ for any 
 code of Type $\rho $. 
In particular the
 full weight enumerators of codes of Type $\rho $ are invariant under
${\mathcal C}(\rho )$ (see \cite[Theorem 5.5.1]{NRS}). 

\begin{remark}\label{WEC}
The {\em Weight Enumerator Conjecture} \cite[Conjecture 5.5.2]{NRS} 
states that in this general situation the fixed space of ${\mathcal C}(\rho )$
is spanned by the full weight enumerators
of codes of Type $\rho $. 
\end{remark}

In fact we do not know a counterexample and 
\cite[Theorem 5.5.5 and Theorem 5.5.7]{NRS} assert the 
truth of the Weight Enumerator Conjecture
for fairly large classes of finite form rings 
including matrix rings over finite fields.

We recall the action of the associated Clifford-Weil 
group for the representation 
$\rho ^{{\bf a}}$ from Remark \ref{sum}
with respect to the $\C $-basis 
$(b_w : w = (w_1,\ldots, w_N) \in V^N)$ 
of the group algebra of $V^N$:
 $${\mathcal C} (\rho^{{\bf a}} ) = \langle
m_r, d_{\phi }, h_{{\iota},u_{\iota},v_{\iota} } : r\in R^{\times }, \phi \in \Phi , {\iota}
=u_{\iota}v_{\iota} \mbox{ sym. idem. }  \rangle $$
where
$$
m_{r}: b_{w} \mapsto b_{rw},\hspace{0.5cm}
d_{\phi}: b_{w} \mapsto \prod _{j=1}^N \exp(2 \pi i \rho_{\Phi }(\phi ) (w_j) )^{a_j} b_{w}
$$
and
$$
h_{{\iota},u_{{\iota}},v_{{\iota}}}: 
b_{w} \mapsto \frac{1}{|{\iota}V|^{N/2}} \sum_{u \in {\iota}V}
\prod _{j=1}^N \exp(2 \pi i\beta (u_j,v_{{\iota}}w_j))^{a_j} b_{u+(1-{\iota})w}.
$$
Here the $\iota $ runs through the set of all $R^{\times }$-conjugacy classes 
of symmetric idempotents in $R$ and $u_{\iota }, v_{\iota } \in R$ are
as in Definition \ref{symid}. 

So the transformations 
$m_r$ are represented as permutation matrices on the chosen basis and the 
transformations $d_{\phi }$ as diagonal matrices. 
It is clear that $\fwe(C)$ is invariant
under all $d_{\phi }$  and all $m_r$, 
if $C$ is an isotropic code in $\rho^{{\bf a}} $. 
For self-dual codes the invariance under $h_{{\iota},u_{{\iota}},v_{{\iota}}} $
follows from a general MacWilliams theorem, see \cite[Example 2.2.6]{NRS}.

\begin{definition} \label{Galois}
	Let $f:=f(\rho )$ and let $\zeta_f:=\exp(\frac{2 \pi i}{f}) \in \C $.
	Put 
	$F_1:= \Q( \{ \sqrt{|\iota V|} : \iota \in R \mbox{ sym. idem.} \} ) $
	and 
	$K(\rho ) := \Q(\zeta _f) F_1$
 the abelian number field containing all entries 
of ${\mathcal C}(\rho )$. 
We choose a complement $F$ of $F_1\cap \Q(\zeta_f) $ in $F_1$. 
For 
	$a \in (\Z/f\Z)^{\times } $  put $\gamma_a$ to 
	denote the Galois automorphism 
of $ K(\rho )$ that is the identity on $ F $
and raises $\zeta _f$ to the $a$th power. 
Then $\gamma : (\Z/f\Z)^{\times } \to \Gal (K(\rho) / F), a\mapsto \gamma_a$ 
is a group isomorphism. 
\end{definition}

\begin{remark}
	Let $a\in (\Z/f\Z)^{\times } $. For a 
	symmetric idempotent $\iota \in R$ we put 
	$\epsilon(\iota, a) \in \{ 1,-1 \} $ such that 
	$\gamma_a(\sqrt{|\iota V|}) = \epsilon(\iota, a) \sqrt{|\iota V|} $.
	Then 
	$${\mathcal C}(\rho ^{(a)}) = \langle 
	m_r, \gamma_a(d_\phi ), \epsilon(\iota,a) \gamma_a(h_{\iota ,u_{\iota}, v_{\iota }} )
	\mid r\in R^{\times }, \phi \in \Phi, 
	\iota \mbox{ sym. idem. } \rangle .$$
\end{remark}

By \cite[Theorem 5.5.3]{NRS} the order of the
scalar subgroup in ${\mathcal C}(\rho )$ is the greatest common divisor of the
lengths of self-dual isotropic codes in $\rho $. 
Combining this with Remark \ref{iota} we obtain 
the following lemma. 

\begin{lemma}
Assume that there is a symmetric idempotent $\iota \in R$ for which 
$|\iota V| $ is not a square in $\Q $.  
	Then $-I_{v} \in {\mathcal C}(\rho )$.
\end{lemma}

So for all $a\in \Z$ prime to $f$ the groups 
${\mathcal C}(\rho ^{(a)}) $ and $ \gamma_a({\mathcal C}(\rho )) $
are equal. 

 \section{Conjugate weight enumerators and the main theorem} 

 We keep the notation of the previous section, in particular 
 $\rho $ is a finite representation of a form ring $(R,\Phi )$ and 
 $K:=K(\rho )$ is an abelian extension of the field $F$ from 
 Definition \ref{Galois}. 
 Put $f:=f(\rho )$ and
 $$\Gamma := \Gal (K/F) =: \{ \alpha _1,\ldots , \alpha _n \} ,$$
 a finite group of order $n$ isomorphic to $(\Z/f\Z)^{\times }$.
 We always assume that there is some length $N_0$ such that 
 there is a self-dual isotropic code in $\rho ^{N_0}$. 
 Then ${\mathcal C}(\rho )$ is a finite subgroup of $\GL_v(K)$,
 where $v=|V|$.
 We also assume that the Weight Enumerator Conjecture 
 \ref{WEC} holds for 
 the form ring $(R,\Phi )$.
For ${\bf a} =(a_1,\ldots , a_N)$
as in Remark \ref{sum} we consider the Clifford-Weil groups ${\mathcal C}(\rho ^{{\bf a}})$ 
introduced in Section \ref{CliffordWeil}.

\begin{definition}
	We say that ${\bf a}$ satisfies the {\em sign condition}
	if $$\prod _{i=1}^N \epsilon (\iota , a_i) = 1 $$ for all 
	symmetric idempotents $\iota \in R$.
\end{definition}

For $1\leq j \leq n$ we put 
$$ D_j := \{ i \in \{ 1,\ldots , N\}  \mid  \gamma _{a_i} = \alpha _j \} \mbox{ and } d_j:=|D_j| \in \N _0.$$ 
Put ${\bf d} := (d_1,\ldots , d_n)$.

\begin{definition}\label{ccweDef}
	The {\em $\Gamma $-conjugate complete weight enumerator} 
	$$\ccwe(C) := \ccwe^{{\bf a}}(C) \in 
	K[ {\bf x} \circ \Gamma ]_{{\bf d} } $$
	of a code $C$ in $\rho ^{{\bf a}}$ 
	is defined as 
	$$\ccwe(C) := \sum _{c\in C} \prod _{j=1}^n \prod _{i\in D_j} x_{c_i}\circ \alpha_j .$$
	Two codes $C,D$ in $\rho ^{{\bf a}}$ are called 
	{\em equivalent} 
	if there are permutations $\pi_j $ of $D_j$ ($1\leq j \leq n$) 
	such that 
	$D=C \circ \pi $ where $\pi = \pi_1\times \ldots \times \pi_n$.
\end{definition}

\noindent
We are now in the position to state and prove the main result of 
this paper. 

 \begin{theorem} \label{main}
	 Assume that $\rho $ satisfies the Weight Enumerator Conjecture 
	 \ref{WEC} and that 
	 ${\bf a}$ satisfies the 
	 sign condition. Then 
	 the space of invariants of ${\mathcal C}(\rho ) $ in 
	 $K[{\bf x} \circ \Gamma ]_{{\bf d} } $ 
 is spanned by the $\Gamma $-conjugate complete weight enumerators of 
	 self-dual isotropic codes in $\rho ^{{\bf a}}$. 
 \end{theorem} 

 \begin{proof} 
	 Define a $K$-linear map 
	 $$\sigma : K V^N \to K[{\bf x} \circ \Gamma ]_{{\bf d }} , (v_1,\ldots , v_N) \mapsto \prod _{j=1}^n \prod _{i\in D_j} x_{v_i}\circ \alpha_j .$$
If $\prod_{i=1}^N \epsilon (\iota, a_i) = 1$ for
         all symmetric idempotents $\iota \in R$ then 
	 $\sigma $ is a ${\mathcal C}(\rho )$-module epimorphism. 
	 As ${\mathcal C}(\rho )$ is a finite group and 
	 $\Char (K) = 0$, both modules are semisimple. 
	 In particular the space of invariants of ${\mathcal C}(\rho ) $ in
	 $K[{\bf x} \circ \Gamma ]_{{\bf d }} $
	 is the image under $\sigma $ of the ${\mathcal C}(\rho)$-fixed space
	 in $K V^{N} $. 
	 By the Weight Enumerator Conjecture \ref{WEC} 
	 this fixed space is spanned by the full weight enumerators 
	 of self-dual isotropic codes $C$ in $\rho ^{{\bf a}}$.
	 Clearly $\sigma (\fwe(C)) = \ccwe (C) $, so the 
	 $\Gamma $-conjugate complete weight enumerators of these codes span the
	 space 
	 $K[{\bf x} \circ \Gamma ]^{{\mathcal C}(\rho  )}_{{\bf d} } $
	 of homogeneous $\Gamma $-conjugate ${\mathcal C}(\rho )$-invariant polynomials 
	 of degree ${\bf d}$.
 \end{proof}

 \begin{remark}
	 The sign condition in Theorem \ref{main} 
is necessary, as the following easy example shows: 
	 Take $R=\F_5$, the field with 5 elements. 
	 Then $J=\id $. We take $\Phi = \leftBra M \rightBra $.
	 We specify the representation $\rho $ by putting $V:=\F_5 = \Z/5\Z$ and 
	 $\beta(x,y) := xy/5 \in \Q/\Z$. Then $K(\rho ) = \Q(\zeta _5)$ and 
	 $F=\Q $.
	 Consider the representation 
	 $\rho ^{(1,2,2,2)}$. As $\epsilon (1,1) = 1$ and $\epsilon (1,2) = -1$,  the sign condition  is not satisfied. 
	 Put 
	 $\Sigma := 4b_{(0,0,0,0)} - \sum _{v\in I} b_v $, where 
	 $$I:= \{ 0\neq (v_1,v_2,v_3,v_4) \in \F_5^4 \mid v_1^2+2(v_2^2+v_3^2+v_4^2) = 0 \} $$ 
	 is the set of isotropic vectors in $\rho ^{(1,2,2,2)}$.
	 Then $m_r(\Sigma ) = d_{\phi }(\Sigma ) =\Sigma $ for 
	 all $r\in \F_5^{\times }$ and $\phi \in \Phi $ but 
	 $h_{1,1,1}(\Sigma ) = - \Sigma $. 
	 The map $\sigma $ from the proof of Theorem \ref{main} 
	 hence maps $\Sigma $ to $0\neq \sigma (\Sigma ) \in K[{\bf x} \circ \Gamma ]^{{\mathcal C}(\rho )} _{1,3,0,0} $. 
	 This invariant space is of dimension 1.
	 However, as the discriminant of the 
	 quadratic form $x_1^2+2(x_2^2+x_3^2+x_4^2) $ is not a
	 square in $\F_5$,  there are 
	 no self-dual isotropic codes in $\rho ^{(1,2,2,2)}$. 
 \end{remark}

 The same method gives $h_{1,1,1}$-anti-invariants for primes $p\equiv _4 1$.

 \subsection{Higher genus Clifford-Weil groups} 

 The  Clifford-Weil groups associated to a representation $\rho $ 
 form an infinite sequence of matrix groups ${\mathcal C}_m(\rho )$ in exponentially 
 growing dimension $v^m$ for $m\in \N$ that comes with  ring epimorphisms 
 $\Phi _m: \Inv ({\mathcal C}_m(\rho )) \to \Inv ({\mathcal C}_{m-1}(\rho ))$  
 that are isomorphisms on the spaces of invariants of small degree. 
 A similar phenomenon occurs for  the space of $\Gamma $-conjugate invariants.

 The full genus-$m$ weight enumerator $\fwe_m(C)$ of $C$ in $\rho $ is the 
 sum over all $m$-tuples of 
 code words of $C$. For $C\leq V^N$ we can think of such an $m$-tuple 
 $c\in C^m$
 as a matrix $c \in V^{m\times N}$,  where the rows correspond to 
 code words. For $1\leq i \leq N$ the $i$th column $c_i$ of this 
 matrix is an element of $V^m = R^m\otimes V$.
In the notation of Definition \ref{ccweDef} we denote by
        $$\ccwe_m(C) := \sum _{c\in C^m} \prod _{j=1}^n \prod _{i\in D_j} x_{c_i}\circ \alpha_j $$ 
 the
        {\em  $\Gamma $-conjugate complete weight enumerator of genus $m$}
	of the 
	code $C$ in $\rho ^{{\bf a}}$.
Now the alphabet of $C^m$ is $V^m$ and hence 
 a module over the ring $\Mat _m(R)$ of $m\times m$-matrices over
 $R$. This gives rise to a notion of a matrix ring of a form ring
 (see \cite[Section 1.10]{NRS}) such that $R^m\otimes \rho $ is 
 a representation of this matrix ring. 
 The {\em genus-$m$ Clifford-Weil group}  is 
 $${\mathcal C}_m(\rho ) := {\mathcal C}(R^m \otimes \rho ) \mbox{ (see \cite[Definition 5.3.4]{NRS}).} $$
 Denote by $K[{\bf x^{(m)}} ]$ the symmetric algebra of $V^m$. 
 Theorem \ref{main} has the following immediate consequence. 

 \begin{corollary}
	 Assume that $\rho $ satisfies the Weight Enumerator Conjecture 
	 \ref{WEC} and that 
          ${{\bf a}}$ satisfies the
         sign condition. Then
	 the space of invariants of ${\mathcal C}_m(\rho  ) $ in
	 $K[{\bf x^{(m)}} \circ \Gamma ]_{{\bf d}} $
 is spanned by the $\Gamma $-conjugate complete weight enumerators of genus $m$
	 of
	 codes of Type $\rho  ^{{\bf a}}$.
 \end{corollary}

 \begin{remark}
	 The map $\Phi _m$ mapping $x_{(v_1,\ldots , v_m)^{tr}}$ 
	 to $x_{(v_1,\ldots , v_{m-1})^{tr}  }$ if $v_m=0$ and 
	 to 0 if $v_m\neq 0$ defines a ring epimorphism 
	 $$
	 \Phi_m:K[{\bf x^{(m)}} \circ \Gamma ]
	 \to K[{\bf x^{(m-1)}} \circ \Gamma ] $$ 
	 that preserves the grading from Section \ref{Secconjpol} and maps 
	 $\ccwe _m(C)$ to  $\ccwe_{m-1}(C)$ for any code $C$ in $\rho ^{{\bf a}}$.
	 In particular it induces an epimorphism from 
	 the space of invariants of ${\mathcal C}_m(\rho  ) $ in
	 $K[{\bf x^{(m)}} \circ \Gamma ]_{{\bf d} } $
	 to the space of invariants of ${\mathcal C}_{m-1}(\rho  ) $ in
	 $K[{\bf x^{(m-1)}} \circ \Gamma ]_{{\bf d}} $.
	 \\
	 If all self-dual isotropic codes in $\rho ^{{\bf a}}$ 
	 are generated by $m-1$ elements then this  epimorphism
	 is in fact an isomorphism. 
	 In this case two codes are equivalent (see Definition \ref{ccweDef}) 
	 if and only if they have the same genus-$m$ weight enumerator 
	 and these weight enumerators of the equivalence classes of 
	 codes of Type $\rho ^{{\bf a}}$ form a basis of 
	 the space of invariants of degree ${\bf d}$. 
 \end{remark}

 \section{Complex conjugate invariants} \label{examples}

In this section we give the most important example where  
all $a_i $ from Remark \ref{sum} are  $\pm 1$.
As $\gamma _{-1} $ is the complex conjugation $\overline{\phantom{x}}$ 
and the field 
$F_1$  from Definition \ref{Galois} is totally real, we always have
$\epsilon(\iota, -1) = 1$ for all symmetric idempotents $\iota \in R$ and
hence the sign condition is fulfilled. 

In particular the main result of \cite{Bannai} is an immediate 
consequence of Theorem \ref{main}. 
Whereas \cite{Bannai} only considers  doubly even, self-dual binary codes 
Theorem \ref{main} allows to consider complex conjugate invariants for 
more general Clifford-Weil groups. 

In view of possible  applications to complex projective designs 
we are mainly interested in
codes of Type $\rho ^{(1^N,(-1)^N)}$ which we call Type $(N,N)$ for short. 
As we want to have as few invariants as possible we consider those 
representations $\rho $ where the quadratic group $\rho _{\Phi }(\Phi )$ 
is as large as possible: 
For self-dual codes over fields of odd characteristic we focus on 
self-dual codes that contain the all-ones vector ${\bf 1}=(1,\ldots , 1)$ and
in even characteristic on the generalized doubly even codes.

\subsection{Codes of Type $(N,N)$} 

\subsubsection{The trivial codes} 
The standard norm $z_m:=||\ ||^2$ on $\C^{v^m}$ 
is an invariant for the full unitary group 
$\mathrm{U}_{v^m}(\C)$. 
In fact $z_m^N$ is  the genus-$m$ complex conjugate 
weight enumerator of the code $T_{N,N}$ 
with generator matrix $(I_N|I_N)$,  the trivial code of Type $(N,N)$.

\subsubsection{Decomposable codes}
For codes $C_i$ of Type $(N_i,N_i)$, $i\in \{ 1, 2\}$, 
the orthogonal sum 
$C_1\oplus C_2$ is a code of Type $(N_1+N_2,N_1+N_2)$. 
If $(A_i|B_i)$ is a generator matrix for $C_i$ then 
$C_1\oplus C_2$ has generator matrix 
$\begin{pmatrix} A_1 & 0 & B_1 & 0 \\ 0 & A_2 & 0 & B_2 \end{pmatrix} $. 
Clearly the genus-$m$ weight enumerator of $C_1\oplus C_2$ is 
just the product of the genus-$m$ weight enumerators of $C_1$ and $C_2$.

In the following we only list representatives of the
equivalence classes of  {\em indecomposable} codes,
i.e.,  codes that cannot be written as a non-trivial orthogonal sum.

\subsubsection{The doubling construction}
There is another general construction of codes in $\rho ^{(N,N)}$ which 
we call the {\em doubling construction}:
\\
Let $C$ be a self-orthogonal isotropic code $C$ in $\rho ^{N}$
and put $D:=C^{\perp} \supseteq C $. Then
$$\Double(C):= \{ (d,d+c)  \mid d\in D, c\in C \} $$ 
is a code of Type $\rho^{(N,N)}$.
If $C$ is a self-dual isotropic code in $\rho ^{N}$, i.e. $C=D$, then
$\Double(C) =\{ (c_1,c_2) \mid c_1,c_2 \in C \}$.

The trivial code $T_{N,N}$  is the double of the zero code,
$T_{N,N} = \Double (\langle 0^N \rangle )$.

\subsection{Enumeration of equivalence classes of codes} 

In this section we sketch some methods to enumerate all 
self-dual isotropic codes in $V^N$. 

\subsubsection{The mass formula} \label{mass}
In many situations the unitary group (in the examples below 
the orthogonal group of the quadratic space ${\bf 1}^{\perp }/\langle {\bf 1} \rangle $) acts transitively on the set of codes of Type $(N,N)$;
their number $t_N$ is the number of isotropic subspaces in a certain 
finite geometry.
If $C_1,\ldots , C_{h_N}$ represent the $\mathrm{Sym}_N \times \mathrm{Sym}_N $-orbits of 
codes of Type $(N,N)$ then we obtain the following mass formula:
$$ \frac{t_N}{N! \cdot N!} = \sum _{i=1}^{h_N} \frac{1}{|\Aut(C_i)|} $$ 
where $\Aut (C_i)$ is the stabiliser in $\mathrm{Sym}_N\times \mathrm{Sym}_N$
of the 
code $C_i\leq V^{N+N}$.
This formula can be used to check the completeness of a list of 
pairwise inequivalent codes but also yields a method to enumerate 
these codes by partitioning the orbit of the unitary group into 
$\mathrm{Sym}_N\times \mathrm{Sym}_N$-orbits. 

\subsubsection{Kneser neighbors} 
A more efficient method to enumerate self-dual codes is the 
Kneser neighbor method \cite{Kneser} (see also \cite{Hecke} for an
application to codes and \cite{Voight} for a survey). 

Starting from one self-dual code $C$ (e.g. $C=T_{N,N}$) 
one enumerates all $\Aut(C)$-orbits 
of neighbors, i.e. those equivalence classes of 
codes $D$ of Type $\rho $ that intersect 
$C$ in a maximal subcode. 
Continuing with the neighbors one successively enumerates all 
codes up to equivalence. 

\subsubsection{Using projections} 

Denote by $\pi_1:V^{N+N} \to V^N $ 
the projection onto the first 
$N$ components and by $\pi _2$ the projection onto the last $N$ coordinates.
For a self-dual isotropic code $C$ of Type $(N,N)$ 
the kernel of $\pi _2$ is the dual of the image of $\pi _1$. 

There is a shortcut to classify the codes $C\leq \F_q^{N+N}$ of Type $(N,N)$ 
for which $\pi _1(C) = \F_q^N $. If the involution is trivial
then each such code has a generator matrix 
$(I_N | A )$ with $A A^{tr} = I_N$. Equivalence of codes translates into 
the action of ${\mathrm{Sym}}_N \times {\mathrm{Sym}}_N$ on 
$\mathrm{O}_N (\F_q) := \{ A \in \F_q^{N\times N} \mid A A^{tr} = I_N \} $ 
by $A \cdot (g,h) := g^{-1} A h $.
Representatives of the $\mathrm{Sym}_N$-double cosets in $\mathrm{O}_N(\F_q)$ hence 
yield generator matrices for representatives of the equivalence classes of codes. 
If one only wants to classify codes that contain the all-ones vector of length $2N$, 
one needs to additionally assume that sum of the rows of $A$ is the all-ones vector 
${\bf 1}$ of length $N$, i.e. ${\bf 1} A = {\bf 1}$. 
As ${\bf 1}g = {\bf 1}$ for all permutation matrices $g\in \mathrm{Sym}_N$ this property is
invariant on the $\mathrm{Sym}_N$-double cosets.

\subsection{Doubly even binary self-dual codes} 
This is the most important example and has been considered in \cite{Bannai}. 
The Type of these codes is denoted by $2_{\II}$ in \cite[Section 2.3.2]{NRS},
has conductor $8$, 
and the associated Clifford-Weil group is the complex Clifford group.
The indecomposable codes of length $(N,N)$, i.e. of Type $2 _{\II}^{(N,N)} $ for 
$N=1,2,3,4,5$ are 
$T_{1,1}$ and $ \Double(\langle(1,1,1,1)\rangle)$ which is called
$g_{4,4}$ in \cite{Bannai}. 

\subsection{Doubly even Euclidean self-dual codes over $\F_4$} 

Doubly-even euclidean self-dual codes over fields of characteristic 2 have
been studied in \cite{NQR}.  Their Type $q_{\II}^E$ is discussed in \cite[Theorem 2.3.2]{NRS}. For $q=4$ the conductor is $4$.
Enumerating the equivalence classes of 
codes in $(4_{\II}^E)^{N,N}$ for $N=1,2,3,4$ we 
obtain the following indecomposable ones:

\begin{tabular}{|c|c|c|} 
	\hline
	Length & number & indecomposables \\ 
	\hline
	(1,1) & 1 & $T_{1,1}$  \\
	(2,2) & 1 & \\
	(3,3) & 2 & $\Double(\langle (1,\omega , \omega ^2  ) \rangle ) $ \\
	(4,4) & 5 & $ \Double(\langle(1,1,1,1)\rangle)$, $\Double (Q_4)$, $C_{4,4}(\F_4 )$ \\
	\hline 
\end{tabular} 
\\
Here $Q_4 $ has
generator matrix $\left(\begin{array}{cccc} 1 & 0 & \omega & \omega^2 \\ 
0 & 1 & \omega ^2 & \omega \end{array} \right) $ and 
the generator matrix of $C_{4,4}$ is
$$\left(\begin{array}{cccccccc} 
	1 & 0 & 0 & 0 & 1 & 1 & \omega & \omega ^2 \\ 
	0 & 1 & 0 & 0 & 1 & 1 & \omega^2 & \omega \\
	0 & 0 & 1 & 0 &  \omega^2 & \omega & 1  & 1\\ 
	0 & 0 & 0 & 1 &  \omega & \omega^2 & 1 & 1 
\end{array} \right) .$$

\subsection{Ternary codes that contain the all-ones vector} 

The Type of self-dual ternary codes that contain the all-ones vector is 
$3_1^{E} $, has conductor 3,  and is described in \cite[Section 7.4.1]{NRS}. We obtain the following list of codes for small lengths:

\begin{tabular}{|c|c|c|}
        \hline
        Length & number & indecomposables \\
        \hline
        (1,1) & 1 & $T_{1,1}$  \\
        (2,2) & 1 & \\
	(3,3) & 2 & $\Double(\langle (1,1 ,1 ) \rangle ) $ \\
	(4,4) & 3 & $ \langle I_4 | 2J_4-I_4 \rangle $ \\
        \hline
\end{tabular}
\\
where the last code has generator matrix $(I_4|2J_4-I_4)$ 
and $J_4$ is the all-ones matrix. 

\subsection{Quinary codes that contain the all-ones vector} 

We consider the analogous Type as in the previous section, where 
the underlying field is $\F_5$ and the conductor is 5.
We obtain the following list for small lengths:

\begin{tabular}{|c|c|c|}
        \hline
        Length & number & indecomposables \\
        \hline
        (1,1) & 1 & $T_{1,1}$  \\
        (2,2) & 1 & \\
	(3,3) & 2 & $ \langle I_3 | 4J_3-I_3 \rangle $ \\
	(4,4) & 5 & $ \langle I_4 | 3J_4-I_4 \rangle $, 
	$\Double(\langle (1,2,3,4 ) \rangle ) $, $C_{4,4}(\F_5)$ \\
        \hline
\end{tabular}
\\
where $C_{4,4}(\F_5)$ has generator matrix 
$$\left(\begin{array}{cccccccc} 
	1 & 0 & 0 & 0 & 0 & 3 & 4 & 4 \\ 
	0 & 1 & 0 & 0 & 3 & 2 & 3 & 3 \\
	0 & 0 & 1 & 0 & 4 & 3 & 0 & 4 \\
	0 & 0 & 0 & 1 & 4 & 3 & 4 & 0
\end{array} \right) .$$


%
%
%
%
%
%
%

\section{A Schur-Weyl duality for Clifford-Weil groups} \label{SchurWeyl}

Let $\rho $ be a finite representation of a form ring $(R,\Phi)$ with underlying 
left $R$-module $V$.
Then the associated Clifford-Weil group $G:={\mathcal C}(\rho )$ 
acts on $\C  V $ and by diagonal action on all 
tensor powers $W_N:=\C V^N  = \otimes ^N \C V  $. 
Let $W_N^* := \Hom (W_N,\C) $ denote the dual space of $W_N$. 
Then $W_N^*$ is also a $\C G$-module, where the action is given by 
$$ f\cdot g  : w\mapsto f(w g^{-1}) \mbox{ for all } w\in W_N, g\in G, f\in W_N^* .$$

It is well known that the linear map 
$$ \varphi : W_N\otimes W_N^* \to \End_{\C} (W_N), \mbox{ defined by } (w,f) \mapsto (x\mapsto f(x) w ) $$
for $w\in W_N$ and $f\in W_N^* = \Hom (W_N,\C) $ is an isomorphism of vector spaces.
Then $\varphi $ is an isomorphism of $\C G$-modules where 
$G$ acts on $\End_{\C }(W_N)$ by conjugation. 
The commuting algebra $\End_{\C G}(W_N) $ 
is the fixed space of the $G$-action on $\End_{\C}(W_N) $.
If the Weight Enumerator Conjecture holds for $\rho $ then 
by Theorem \ref{main} this fixed 
space of $G={\mathcal C}(\rho )$ is spanned by the images of the full weight enumerators of self-dual isotropic codes
in $\rho ^{(N,N)}$, so we have the following theorem.

\begin{theorem}\label{schur}
Assume that the Weight Enumerator Conjecture holds for $\rho $.
	Let $C_1,\ldots , C_{t_N} $  be a complete list of 
	self-dual isotropic codes in $\rho ^{(N,N)}$. 
	Then 
	$$(\varphi (\fwe^{(N,N)}(C_1)) , \ldots , 
	\varphi (\fwe^{(N,N)}(C_{t_N})) )$$ is a generating set of 
	$\End_{{\mathcal C}(\rho )}(\C V^N ) $
\end{theorem} 

Let $g_N(\rho )\in \N$ be such that any 
 code of Type $(N,N)$ can be generated by $g_N(\rho )$ elements.
Whenever $m\geq g_N(\rho )$ the genus-$m$ full weight enumerators of codes of length $(N,N)$ are linearly independent.

\begin{corollary} (\cite[Theorem 4.3]{Gross})
	It $m \geq g_N(\rho )$ then the set 
	$$(\varphi (\fwe_m^{(N,N)}(C_1)) , \ldots , 
	\varphi (\fwe_m^{(N,N)}(C_{t_N})) )$$ is a basis of
	$\End_{{\mathcal C}_m(\rho )}(\C (V^m)^N ) $.
\end{corollary}

The number $t_N$ is often well understood (see Section \ref{mass}) and hence 
determines the dimension of the commuting algebra 
of the $N$-fold tensor representations of the genus-$m$ Clifford-Weil groups.
Note that these dimensions are independent of $m$ for $m\geq g_N(\rho )$.

%
%

\end{document}